\documentstyle[12pt]{article}

\newtheorem{thm}{Theorem}
\newtheorem{defn}{Definition:}
\newtheorem{lemma}{Lemma:}

\input{psfig}

\title{The Structure and Enumeration of Link Projections}
\author{Martin Bridgeman$^1$}

\date{October 15, 1994}

\def\pf{{\bf Proof :}}
\def\eproof{$\Box$}

\def\ator{atoroidal }
\def\atorg{atoroidal graph }
\def\atorgs{atoroidal graphs }

\def\nbd{neighborhood }
\def\a{\alpha}
\def\b{\beta}
\def\c{\gamma}
\def\d{\delta}

\def\bdy{boundary }

\begin{document}

\maketitle
\footnotetext[1]{Research at MSRI is supported in part by NSF 
                 grant no.DMS-9022140.}

\begin{abstract}
We define a decomposition of link projections whose pieces we
call atoroidal graphs. We describe a surgery operation on these graphs and 
show that all atoroidal graphs can be generated by performing surgery
repeatedly on  a family of well known link projections. This gives a method of
enumerating \atorgs and hence link projections by recomposing the pieces of the
decomposition.
\end{abstract}

\section{Introduction} 
The problem of enumeration of knots and links has always interested knot 
theorists. 
In this paper we introduce a method of enumerating link projections by first
decomposing them into pieces called
{\em atoroidal graphs}. We define surgery on these atoroidal graphs and  show 
how  they can be enumerated by performing surgery on a well known
family of link projections. By recomposing these atoroidal graphs we can thus
enumerate link projections. I have included an enumeration of atoroidal graphs
to 12 crossings.

A link projection is given by a 4-valent planar graph $G$. To form a link we 
can replace each vertex of $G$ by a crossing. To enumerate links in this way
we must first enumerate link projections. It was Kirkman's success in
enumerating link projections or {\em polyhedra} as he called 
them(\cite{K1},\cite{K2}) that formed the basis of the knot tables of both 
Tait(\cite{T}) and Little(\cite{L1},\cite{L2}). In \cite{C}, Conway 
introduced a notation which made it possible for 
him to enumerate knots to 11 crossings and links to 10 crossings in a single
afternoon where before it had taken years. In his paper Conway defined a
{\em basic polyhedron} to be a polyhedron with no bigon regions and showed 
that every link is obtained by replacing each vertex of a basic polyhedron by 
a rational tangle. These basic polyhedra are closely related to
the atoroidal graphs defined in this paper and can be enumerated 
using the enumeration of atoroidal graphs described.

The decomposition of a link projection into atoroidal graphs is achieved by 
cutting the projection along certain non-trivial curves. We then define 
surgery on an atoroidal graph giving a new atoroidal graph with one more
vertex. This gives a partial ordering on atoroidal graphs  where 
$G_{1} \prec G_{2}$ if a surgery on $G_1$ results in $G_2$ and we show that a 
graph is initial if and only if it has no vertices of a given type. Using this 
we can list the initial objects and thus enumerate all atoroidal graphs by 
repeatedly performing surgery on these initial objects.

The motivation for the paper comes from orbifold theory and hyperbolic 
geometry but a background in these is not necessary here. For a reference see 
\cite{Th}. For  readers interested, these aspects are laid out
in the section on orbifolds.

I would like to thank Curt McMullen, Joe Christy, Rich Schwartz and 
especially my advisor Bill Thurston.

\section{Decomposition}

\subsection*{Link Projections}
Given a link  $L$ a general position projection of $L$ is a 4-valent graph
$G$ embedded in $S^2$. As we are only considering such graphs we will 
use graph to mean a 4-valent graph embedded in $S^2$.
 
\noindent Let $G$ be a graph with vertex set $V$.
\begin{defn}
An {\em n-curve} of $G$ is a simple closed curve in $S^{2}-V$ intersecting $G$
 n times.
\end{defn}

\begin{figure}
\begin{center}~
\psfig{file=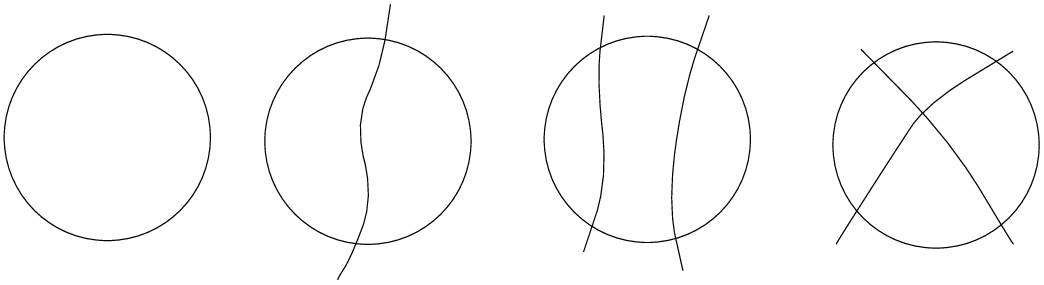,height=0.5in}
\end{center}
\caption{Trivial 0,2,4-curves
\label{trivial024}}
\end{figure}

Let $\a$ be an n-curve($n=0,2,4$) of $G$. Then $\a$ splits $S^2$ into 
two disks 
$D_{1},D_{2}$. We say a component  $D_i$ is {\em trivial} if $D_{i} \cap G$ is
either
empty, a simple arc, two disjoint simple arcs or two arcs crossing at a single
vertex(figure \ref{trivial024}). If $\a$ has a trivial component then $\a$ 
is {\em trivial}. Otherwise $\a$ is called {\em non-trivial}.

\begin{defn}
A graph $G$ is {\em irreducible} all n-curves($n=0,2$) are trivial.
\end{defn}

\begin{defn}
A graph $G$ is {\em atoroidal} if all n-curves($n=0,2,4$) are trivial.
\end{defn}

\subsection*{Decomposition}
If $\a$ is a non-trivial n-curve($n=0,2,4$) of $G$ then we can decompose $G$ 
along $\a$ into  graphs $G_{1}$ and $G_{2}$ as follows. First cut along $\a$, 
this splits the 
sphere  into two disks $D_{1},D_{2}$. To obtain the graph $G_i$ from $D_i$ we 
identify the boundary of $D_i$ to a single point. We say that $G$ decomposes 
into $G_{1}$ and $G_{2}$ 
along $\a$.

We now describe the decomposition of  a link projection $G$ into  
atoroidal graphs.
If  all n-curves($n=0,2,4$) in $G$ are trivial the the decomposition is done. 
Otherwise decompose $G$ into $G_1$ and $G_2$ along a non-trivial 
n-curve($n=0,2,4$) where $n$  is chosen to be as small as possible. Now repeat
the decomposition on the resultant graphs $G_1$ and $G_2$.
It is obvious that this decomposition  terminates.

The decomposition along
non-trivial 0-curves is especially simple, corresponding to splitting a graph 
into its connected components and thus when dealing with connected graphs we
only need concern ourselves with non-trivial 2-curves and 4-curves. 
 
\section{Structure of Atoroidal Graphs}

\subsection*{Almost all Atoroidal Graphs are Hyperbolic}

To investigate the type of \atorgs possible we consider the
cell division of $S^2$ given by $G$. We call the cells,
faces of $G$ and a face $F$ is called an n-gon if it has n vertices of
$G$ on its boundary. Any n-gon $F$ of $G$ gives a 2n-curve $\alpha_F$ by 
taking the boundary
of a small neighborhood $N_F$ of $F$.

\begin{figure}
\begin{center}~
\psfig{file=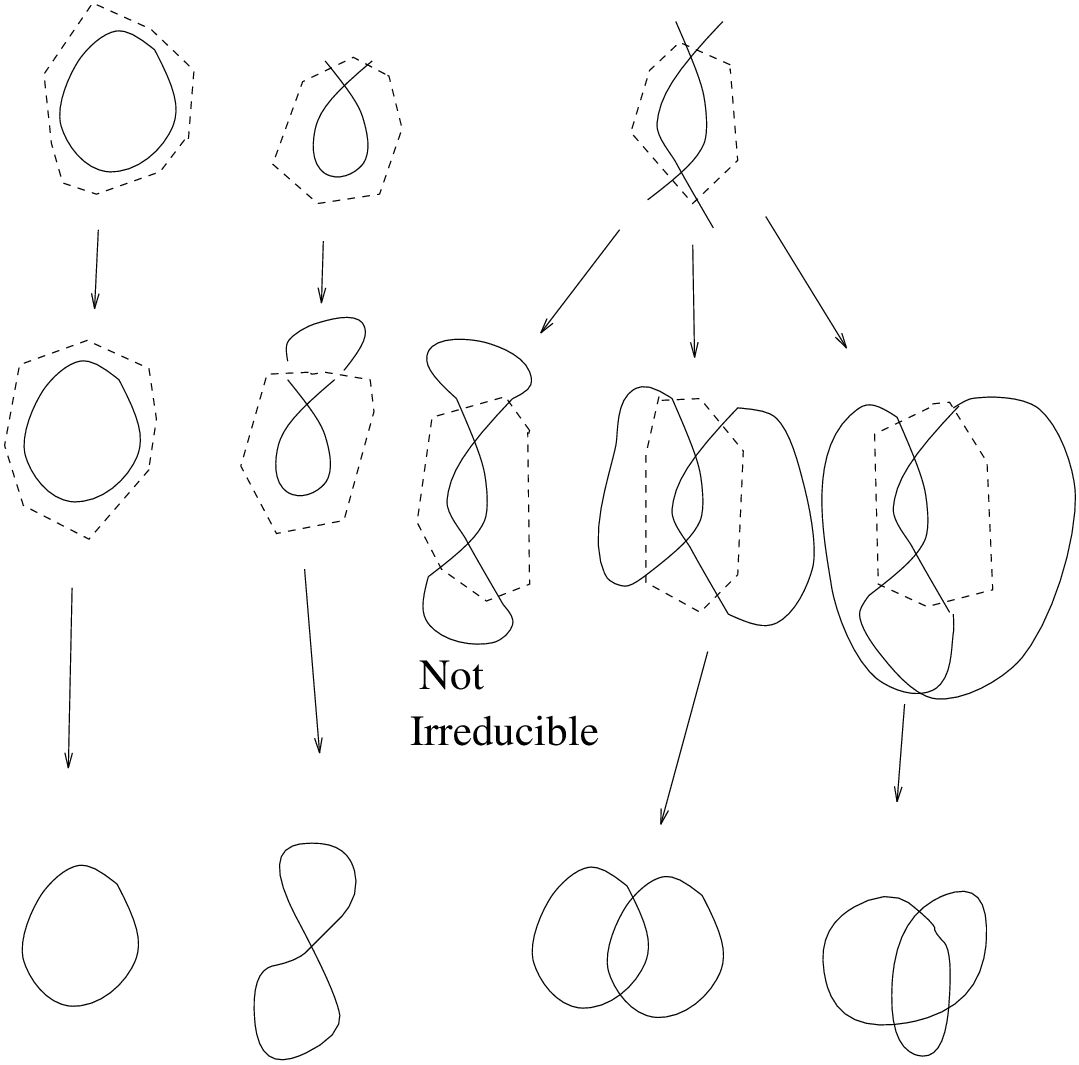,height=3in}
\end{center}
\caption{Exceptions 
\label{exceptions}}
\end{figure}

If $G$ has a 0-gon $F$ then $\alpha_F$  is
a trivial 0-curve (i.e. it bounds a disk $D$ in $S^{2}-G$) as $G$ is
irreducible. Attach
a disk to the boundary of $N_F \Rightarrow  G$ is just
a jordan curve or the unknot projection(figure \ref{exceptions}).

If $G$ has a 1-gon(monogon) $F$ then $\alpha_F$ is a
trivial 2-curve as $G$ is irreducible.
Attach a neighborhood of an arc to $N_F
\Rightarrow G$ is the graph having the form of the number 8  
(figure \ref{exceptions}).

If $G$ has a 2-gon(bigon) $F$ then $\alpha_F$ is a trivial 4-curve as $G$ is 
atoroidal. Attaching a neighborhood of a vertex to $N_F \Rightarrow G$ is the
projection of the trefoil and attaching a neighborhood of two parallel
arcs of $G$ to $N_F$ can be done in 2 ways to give two possible graphs but as
 can be seen in figure \ref{exceptions} only one is also 
irreducible. Thus $G$ is 
either the trefoil projection or the Hopf link projection.

Apart from these 4 exceptional graphs all other \atorgs have faces that are
at least triangular and are called {\em hyperbolic graphs}.

\subsection*{Surgery}
\begin{figure}
\begin{center}~
\psfig{file=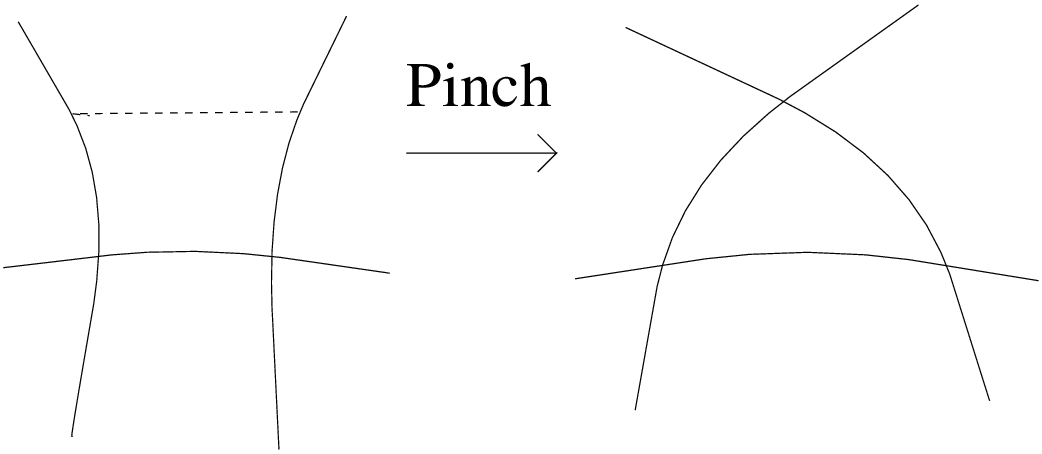,height=1in}
\end{center}
\caption{Surgery
\label{move}}
\end{figure}

Given any \atorg $G$ which has a face $F$ with greater than 3 vertices then we
can perform {\em surgery} on $G$ to give another \atorg $G'$.
We choose edges $e,e_1,e_2$ of $F$ with $e_1,e_2$ adjacent to $e$.
$G'$ is obtained by pinching together $e_1,e_2$, i.e. 
take a simple arc $\alpha \subset F$ with endpoints in the interior
of $e_1,e_2$ resp. and homotope $\alpha$ to a single point. The graph 
$G'$ obtained by performing surgery on $G$ has one more vertex 
than $G$(figure \ref{move}).

\begin{lemma}
$G'$ is atoroidal.
\end{lemma}
\pf Let $\alpha'$ be an n-curve in $G'$ which decomposes $S^2$ into disks
$D'_1,D'_2$. We can get $G$ back by splitting open the new
crossing $v$ and as this splitting can be done in a small
\nbd of $v$ and $\alpha'$ is outside such a neighborhood, we get an
n-curve $\alpha$ in $G$ which splits $S^2$ into disks $D_1,D_2$.
Note that $D'_i$ and $D_i$ are either the same or the former is obtained 
from the latter by pinching two edges together.

As $G$ is \ator then $\alpha$ is trivial for $n=0,2,4$ and
we can assume $D_1$ is trivial.

\noindent $n=0,2 \Rightarrow D'_1 = D_1$ as $D_1$ has at most
one edge intersecting it so we haven't enough edges to pinch 
$\Rightarrow D'_1$ trivial.
$\Rightarrow \alpha'$ trivial.

\noindent $n=4 \Rightarrow D_1$ either \nbd  of vertex or \nbd
of two non-intersecting arcs of $G$. If $D_1$ is \nbd of vertex
  then the only edges
that can be pinched are adjacent $\Rightarrow D'_1 = D_1 \Rightarrow
\alpha'$ trivial. If $D_1$ is \nbd of two parallel arcs
$\Rightarrow D'_1 = D_1$ or $D'_1$ is
obtained by pinching the parallel arcs of $D_1$ together 
$\Rightarrow D'_1$ is \nbd of a vertex $\Rightarrow \alpha'$  
trivial $\Rightarrow G'$ is \ator \eproof \newline

\section{Orbifolds}
This section explains how the work arises out of considering certain orbifolds 
associated with a link projection. In this setting the decomposition and 
surgery we define are the torus decomposition and dehn surgery on these 
orbifolds. An orbifold is a generalization of a manifold in which the space 
is locally modeled on $R^n$ modulo the action of a finite group. For example 
if a group $G$ acts properly discontinuously on a space $M$ then $M/G$ is 
an orbifold and is a manifold if the action is also free. For a reference on 
orbifolds see chapter 13 of \cite{Th}.

\subsection*{Associated Orbifolds}
We associate two orbifolds $O_{G}$ and $O'_{G}$ to a graph $G$ as follows.
We consider $G$ as a graph siting on $S^2$ in $S^{3}$. Let $B$ be a ball
in $S^3$ with boundary $S^2$ and $V$ be the vertex set of $G$.
$O_{G}$ is a polyhedral orbifold with underlying space $X_{O_{G}} = B - V$, 
singular locus $\Sigma_{O_{G}} = S^{2} - V$ and 1-dimensional singular locus
$\Sigma^{1}_{O_{G}} = G - V$. The 1-dimensional singular locus is marked with
$D_{2}$ indicating that any point on it is modeled by
 $D^{3}/D_{2}$ where $D_{2}$
acts by two reflections in planes meeting in right angles.
$O'_{G}$ has underlying space $X_{O_{G}'} = S^{3} - V$ and singular locus
$\Sigma_{O_{G}'} = G - V$. Here the singular locus is 1-dimensional and is
marked with $Z_{2}$ to indicate any point on it is modeled on $D^{3}/Z_{2}$
where $Z_{2}$ acts by rotation of order two.

$O_{G}'$ is the double of $O_{G}$ in the sense of orbifolds.

\subsection*{Decomposition and Surgery}
In \cite{B} we show that the torus decomposition on the orbifolds $O_{G}$
and $O'_{G}$ is the decomposition we've described on $G$.  By Andreev's 
theorem(see \cite{Th}) if $G$ is a hyperbolic graph then $O_G$ can be 
realized as an ideal hyperbolic polyhedron with all dihedral angles right 
angles. Taking the subgroup
of orientation preserving elements of $\pi_{1}(O_{G})$ shows 
us that $O'_G$ can also be realized as a hyperbolic orbifold.

We show(\cite{B}) that for any graph $G$  the double cover of $O'_G$ is a 
link compliment denoted by $L_G$ with one 
component for each vertex of $G$ and is a hyperbolic link complement iff
$G$ is a hyperbolic graph. Also if $G'$ is obtained by surgery on $G$ then 
link complement  $L_{G'}$ is obtained from $L_{G}$ by removing a simple
closed curve, i.e. by {\em dehn drilling}.

\section{Partial Ordering}
Surgery gives \atorgs a partial ordering $\prec$ by defining
$G_1 \prec G_2$ iff $G_2$ is obtained by performing $r$ successive
surgeries on $G_1, r=0,1,2,\ldots$ Note that surgery cannot be performed on
any of the four exceptional \atorgs  and they are never the resultant
graph of surgery, therefore they are isolated objects(both initial and final).
Thus $\prec$ restricts to a p.o. on hyperbolic
graphs.  To study  $\prec$ we show that the initial objects
are a well known family of
graphs and we can generate all \atorgs by performing surgery on these initial 
objects.

\begin{figure}
\begin{center}~
\psfig{file=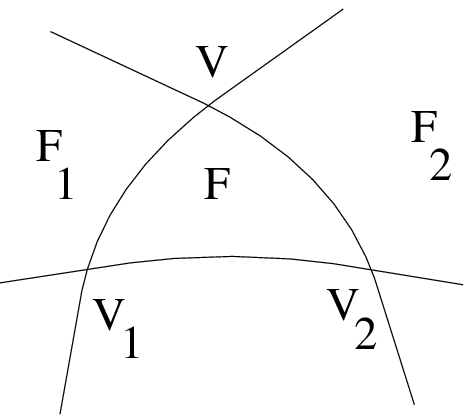,height=1in}
\end{center}
\caption{Simple vertex
\label{local picture}}
\end{figure}

\begin{figure}
\begin{center}~
\psfig{file=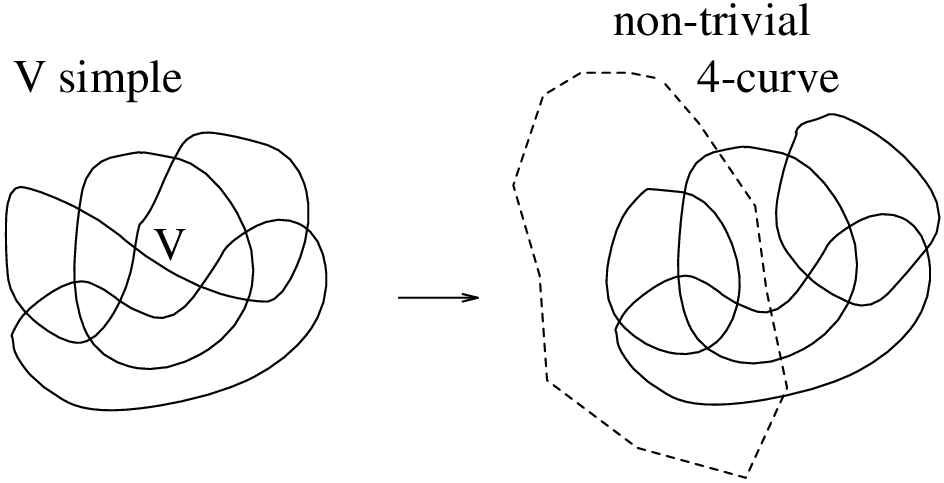,height=1.5in}
\end{center}
\caption{Cutting open at  a simple vertex doesn't work
\label{exception}}
\end{figure}

After surgery has been performed on a graph to give a graph $G$ with a new
vertex $v$. The vertex $v$ is a vertex of a 
triangle $T$ which has adjacent faces $F_1,F_2$ meeting at $v$  each 
being greater than triangular(figure \ref{local picture}). A vertex with this
local structure we call {\em simple}. To find a 
$\bar G$ s.t. $\bar G \prec G$ we might
just look for a simple vertex $v$ and cut open at $v$(there is a
unique way to cut open a  simple vertex) but this doesn't
necessarily give an \atorg as the resultant may have non-trivial
4-curves(figure \ref{exception}). What we will show is that if a graph
has a simple vertex $v$ belonging to a  triangle $T$ then the graph can
be cut open at {\em some} vertex of $T$ to give an atoroidal graph. 
This implies that an initial object cannot have the any simple vertices.
Before proving the stated result we will use it to show what the initial 
objects are.

\begin{figure}
\begin{center}~
\psfig{file=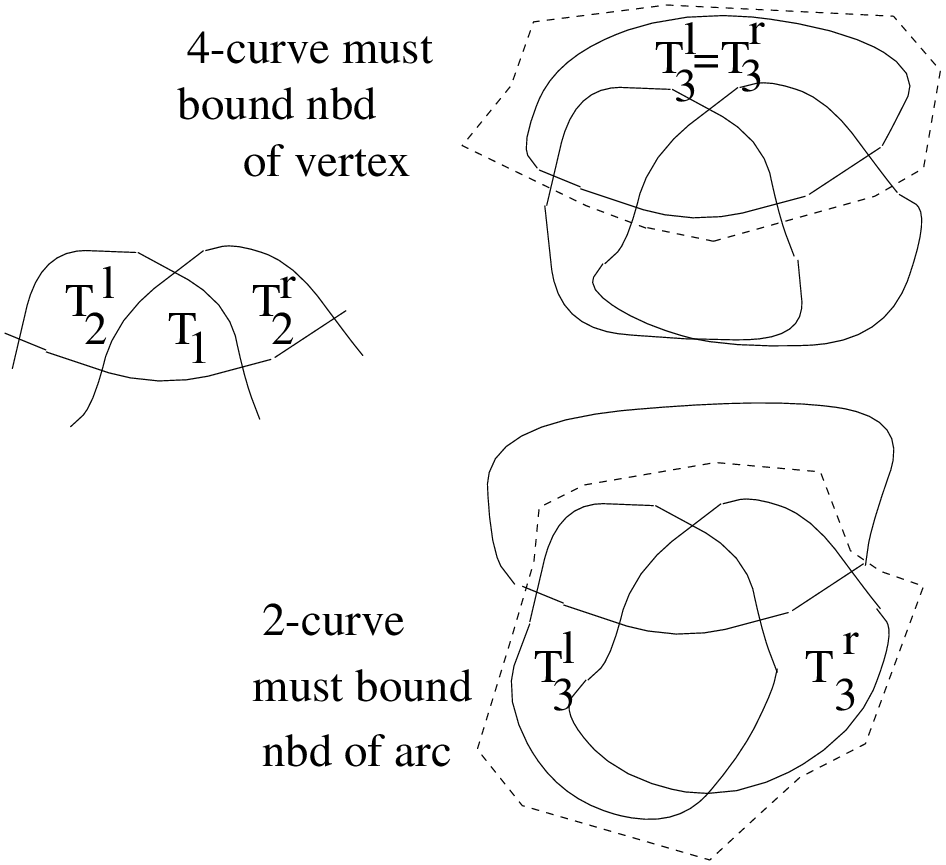,height=1.75in}
\end{center}
\caption{Initial Setup and the resulting  Borromian Rings
\label{borr}}
\end{figure}

Since the exceptions are isolated, they never arise in a sequence
of surgeries and all other initial objects are hyperbolic. Let $G$ be a 
hyperbolic
initial object(not one of the exceptions), calculating the euler number
of the cell division of $S^2$
into faces of $G$ we see that $G$ has a triangular face $T_1$.
As $G$ is initial, $T_1$ has two adjacent triangular faces
$T_2^l,T_2^r$(figure \ref{borr}).
Again using the fact that $G$ is initial we 
have that $T_2^l,T_2^r$ both have a neighboring triangular face
other than $T_1$ labeled $T^l_3,T^r_3$ respectively.

If  $T_3^l = T_3^r$ then $G$  \ator $\Rightarrow G$
is borromian ring projection. Also if $T^l_3,T^r_3$
have a common vertex then $G$ irreducible $\Rightarrow G$
is again the Borromian ring projection which we call $T_3$
(figure \ref{borr}).

\begin{figure}
\begin{center}~
\psfig{file=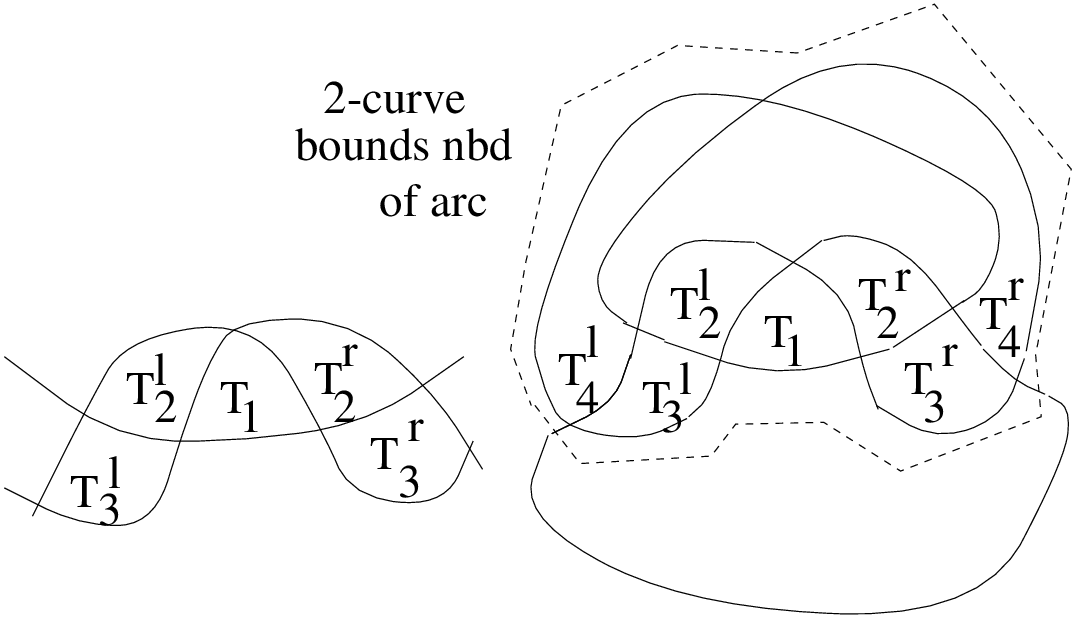,height=1.25in}
\end{center}
\caption{Next Stage and the resulting graph $T_4$
\label{t4}}
\end{figure}

If $T^l_3,T^r_3$ are disjoint then each has a neighboring
triangular face $T^l_4,T^r_4$ other than the previous faces
$T^l_2,T^r_2$. These are unique as the face adjacent to both
$T^l_3$ and $T^r_3$
is at least 4 sided. If now $T^l_4,T^r_4$ have a common vertex
then using $G$ \ator $\Rightarrow G$ is of the form given in 
figure \ref{t4} which we call $T_4$.

\begin{figure}
\begin{center}~
\psfig{file=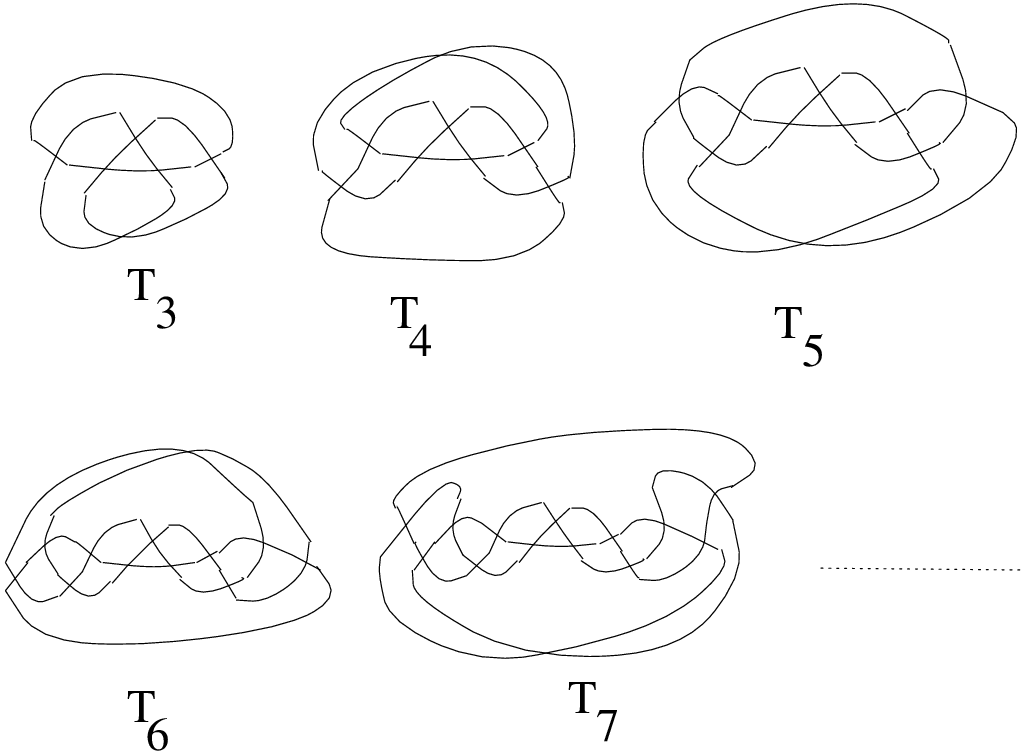,height=2in}
\end{center}
\caption{$T_3$	$T_4$	$T_5$	$T_6$	$T_7$	$\ldots$
\label{tn}}
\end{figure}

Continuing this we get the collection of  
graphs $\{T_n\}_{n \ge 3}$(figure \ref{tn}) which together
with the exceptional \atorgs are
the initial objects of $\prec$ of which only $\{T_n\}_{n \ge 4}$
are  non-terminal. Knowing the initial objects allows us enumerate all \atorgs
by performing surgery repeatedly. 

\begin{lemma}
Let $G$ be an  \atorg  with a simple vertex $v$ of triangle $T$ and  let $G'$ 
be the graph
obtained by cutting $G$ open at $v$. Let $e'_1,e'_2$ be the two edges of face
$F'$ in $G'$ pinched to get $G$  then

\noindent {\bf 1.}  $G'$  is irreducible. 

\noindent {\bf 2.} Any non-trivial
4-curve $\alpha'$ of $G'$ intersects the face $F'$ in a single arc $\beta'$ 
which separates $e'_1,e'_2$.
\end{lemma}

\pf As above let $e'_1,e'_2$ be the two edges of $F'$ pinched together
to give $G$ with $e'$ the edge adjacent to both.
If $\exists$ \nbd $N_{e'}$ of $e'$ s.t.
$N_{e'} \cap \alpha'$ is empty $\Rightarrow$ can pinch $e'_1,e'_2$
in $N_{e'}$ with $\alpha'$ giving  an n-curve $\a$ in $G$.
If $n = 0,2,4$ then $\alpha$ is trivial and splits $S^2$ into
two disks $D_1,D_2$ with $D_1$ trivial. Similarly $\alpha'$ splits
$S^2$ into $D'_1,D'_2$ with either $D'_1 = D_1$ or $D'_1$ obtained from
$D_1$ by cutting open a crossing. In either case this gives $D'_1$
trivial $\Rightarrow$ $\alpha'$ trivial.

\noindent Therefore every non-trivial n-curve ($n= 0,2,4$) must intersect
$F'$ in an arc $\b'$ that has one endpoint on $e'$ and the other
on another edge $e'_3$ of $F'$ with $e'_3 \neq e',e'_1,e'_2$.

If $n = 0$ then $\a'$ doesn't intersect $G' \Rightarrow \a'$ trivial.

\begin{figure}
\begin{center}~
\psfig{file=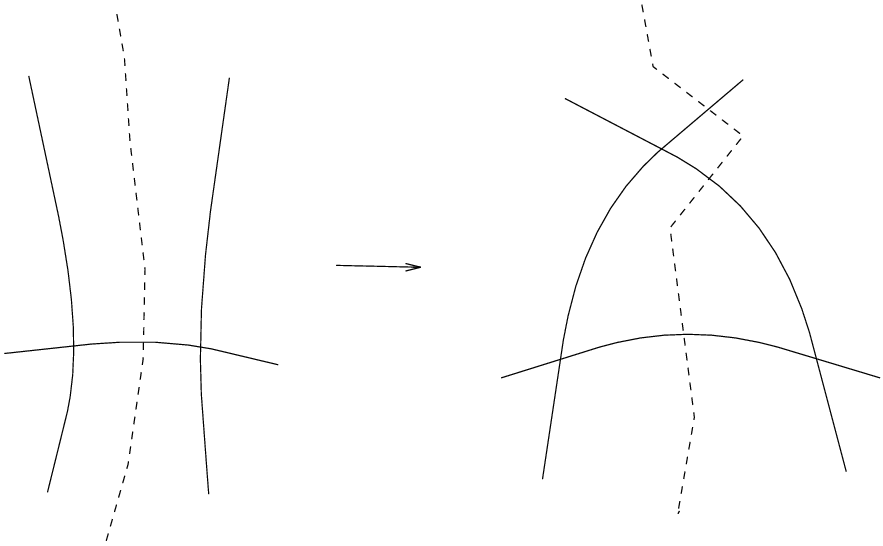,height=1in}
\end{center}
\caption{Pinch through $\a'$
\label{pinch}}
\end{figure}

If $n = 2$ then $\a'$ only intersects $G'$ at the two endpoints of $\b'$.
If we pinch $e'_1,e'_2$ together to get $G$ we can do so by either
 leaving $e'_1$ fixed and pulling $e'_2$ through $\a'$ or vise-versa.
This gives  4-curves $\a^r,\a^l$  resp. in $G$ 
which are identical with $\a'$ 
outside a \nbd of the new vertex $v$ and either go right or left around
$v$ as the curves approach $v$ from inside $T$(figure \ref{pinch}).
 $G$ \ator $\Rightarrow
\a^r$ is trivial $\Rightarrow \a^r$ splits $S^2$ into $D_1,D_2$ s.t.
$D_1$ trivial. The region containing $v$ also contains another vertex
of $T$ so it can't be trivial. Therefore $D_1$ is the region containing 
vertex $v_2$ of $T$ and is \nbd of $v_2 \Rightarrow$  $v$ not simple
as one of faces is bigon. This contradiction
implies that $G'$ has no non-trivial 2-curves $\Rightarrow G'$ is
irreducible.

\begin{figure}
\begin{center}~
\psfig{file=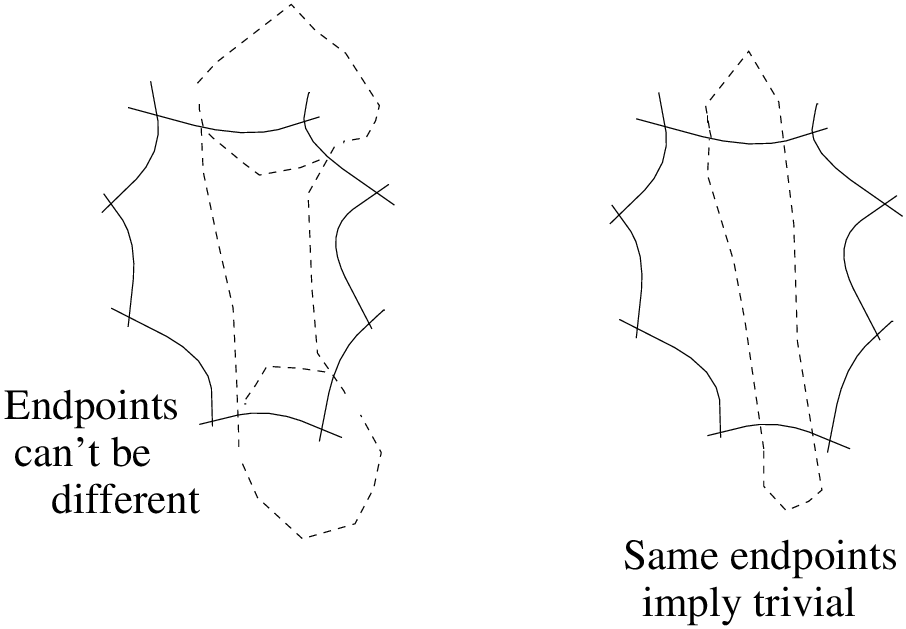,height=1.5in}
\end{center}
\caption{Non-Trivial curve can't intersect in more than one arc
\label{oneint}}
\end{figure}

If $n = 4 \Rightarrow \a' \cap F'$ consists of either 1 or 2 arcs. If it
is 2 arcs $\b'_1,\b'_2$ then traversing around $\a'$ we have 4 connected
arcs $\b'_1,\c'_1,\b'_2,\c'_2$.
We can join endpoints of $\c'_1$ by another arc $\d'_1$ in $F'$ s.t.
$\c'_1 \cup \d'_1$ is a 2-curve in $G' \Rightarrow$ trivial. If endpoints
of $\c'_1$ belong to different edges then both components of 
$S^{2} - \c'_1 \cup \d'_1$ contain vertices which contradicts 
it being trivial. Therefore both
$\c'_1,\c'_2$ are contained in adjacent faces to $F'$ and $\a'$ is
\nbd of two parallel arcs of $G' \Rightarrow \a'$ trivial(figure \ref{oneint}).
Therefore any non-trivial 4-curve in $G'$ intersects $F'$ in a single arc
$\b'$ separating the two edges pinched in $G'$ to obtain $G$ \eproof

\begin{defn}
An n-curve $\a$ ($n \geq 4$) in $G$ is trivial iff either

\noindent $\bullet$ $\a$ is the \bdy of a \nbd of a vertex of $G$

\noindent or  

\noindent $\bullet$ $\exists$ arc $\b$
 intersecting $G$ at most once s.t. $\b \cap \a = \partial{\b}$ and
$\partial{\b}$ splits $\a$ into $\a_1,\a_2$ each containing at least 
two points of $G$. 

\noindent $\b$ is called a compression of $\a$.
\end{defn}

\begin{figure}
\begin{center}~
\psfig{file=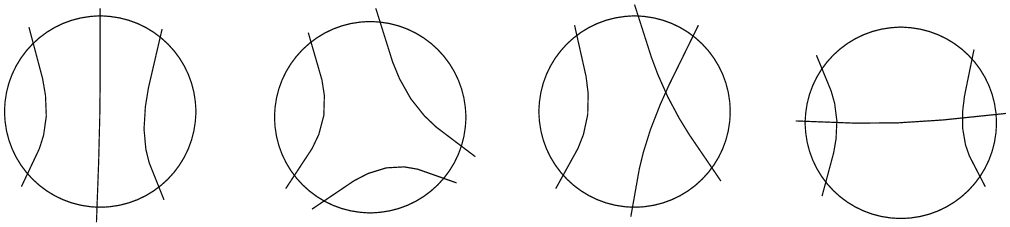,height=0.5in}
\end{center}
\caption{Only trivial 6-curves in \atorg
\label{triv6}}
\end{figure}

If $G$ is \ator then the only trivial 6-curves can easily shown to be
those curves that split $S^2$ into two disks one of which has
one of four types given in figure \ref{triv6}. Disks of this kind
bounding a 6-curve are called trivial.

\begin{lemma}
If $\a'$ is a non-trivial 4-curve in $G'$ then associated with it are two
non trivial 6-curves $\a^l,\a^r$ in $G$. 
\end{lemma}

\pf $\a'$ intersects $F'$ in a single arc $\b'$ separating edges $e'_1,e'_2$
and as before we can pinch $e'_1,e'_2$ together to get graph $G$. This can be
done in two ways, either by fixing $e'_1$ and pushing $e'_2$ across $\b'$
or vise-versa. We get two 6-curves $\a^r,\a^l$ in $G$ both identical to $\a'$
outside a neighborhood of the new vertex $v$ and either going right or left
around the vertex inside the \nbd of $v$ as before(figure \ref{pinch}).

$\a^r$ splits $S^2$ into disks $D_1,D_2$ with $D_1$ containing
vertices $v,v_1$ of $T$ and $D_2$ containing vertex $v_2$. Therefore if
$D_1$ is trivial then it must be the same type as the fourth disk 
described in figure \ref{triv6}. But then $\a'$ would bound a \nbd of
 a vertex of $G$ which contradicts $\a'$ non-trivial. If $D_2$ is
trivial then it is either the same type as the third or fourth disk in
figure \ref{triv6}. If its the third type then as before $\a'$ would
bound a \nbd of a vertex of $G$ which contradicts $\a'$ non-trivial.
If it is the fourth type then adjacent face $F_2$ of $T$ would be a
triangle contradicting  $v$ being simple.
Therefore $\a^r$ is non-trivial and similarly $\a^l$ \eproof

\begin{figure}
\begin{center}~
\psfig{file=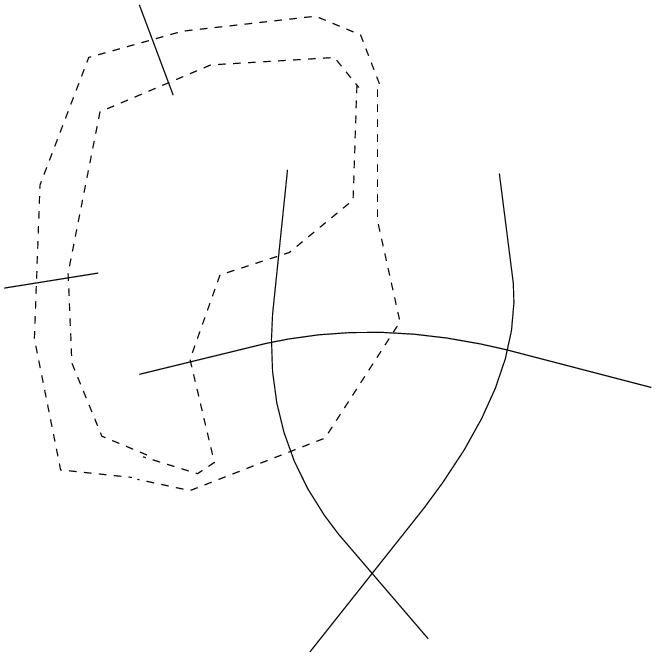,height=1in}
\end{center}
\caption{All vertices of $T$ are simple
\label{nottri}}
\end{figure}

\begin{thm}
If $G$ has a simple vertex $v$ of a triangle $T$ then either

\noindent $\bullet$ splitting open at $v$ gives an \atorg $G'$ 

\noindent or

\noindent $\bullet$  both other vertices $v_1,v_2$ of $T$ are simple and 
splitting at either gives an atoroidal graph.
\end{thm}
\pf If $G'$ is not \ator $\Rightarrow \exists \a'$ non-trivial
4-curve in $G'$ and $\a^r,\a^l$ non-trivial 6-curves in $G$.
Triangle $T$ has adjacent faces $F_1,F_2,F_3$ with both $F_1,F_2$
non-triangular as $v$ is simple. If $F_3$ was triangular
then $\a^r$ splits $F_3$ in two, one piece containing just one vertex
say $v_1$ and the other containing two. Therefore $\a^r$ takes a 
clockwise path about $v_1$ from $T$ through $F_3$. If instead we
take an anticlockwise path $\a'$ gives us another 4-curve in
$G'$ called $\bar{\a}'$(figure \ref{nottri}). Since $\bar
{\a}'$ doesn't intersect $e'$ then it is trivial. Therefore it splits
$S^2$ into disks $\bar{D}'_1,\bar{D}'_2$ with $\bar{D}'_1$
trivial. $\a'$ splits $S^2$ into $D'_1,D'_2$ with $D'_1$ obtained from 
$\bar{D}'_1$ by crossing two adjacent ends so if $\bar{D}'_1$
is \nbd of two parallel arcs then $D'_1$ is \nbd of vertex and if 
$\bar{D}'_1$ is \nbd of vertex then $G'$ would contain a bigon.
Therefore $F_3$ must be non-triangular and both  $v_1$ and $v_2$ are simple.

\begin{figure}
\begin{center}~
\psfig{file=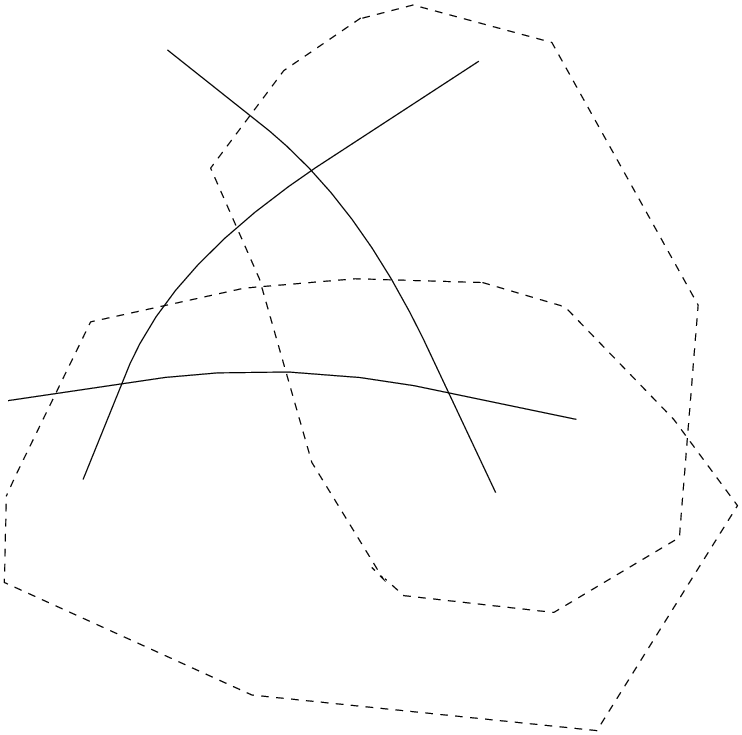,height=1.5in}
\end{center}
\caption{Two intersecting Non-trivial 6-curves
\label{setup}}
\end{figure}

If splitting at $v_1$ doesn't give an \atorg then there is a non-
trivial 4-curve $\a'_1$ in $G'_1$ and non-trivial 6-curves
$\a^r_1$, $\a^l_1$ in $G$. Considering the 6-curves $\a,\a_1$,
where $\a = \a^l$ and $\a_1 = \a^r_1$ we will show that they can be isotoped
to intersect in only two points.Then by showing that they
cannot intersect in the given way(figure \ref{setup}) the theorem
is proven.

Firstly we will show that $\a,\a_1$ can be isotoped to only 
intersect twice. If they intersect any more then $S^{2} - \a \cup \a_1$
contains at least four regions that are disks with boundary
consisting of one arc of $\a$ and $\a_1$. As $\a \cup \a_1$ has 12
intersections with $G$ then each of these four regions cannot
have boundaries being n-curves $n \geq 4$. Therefore one of
these regions $D$ has $\c = \partial{D}$ either a trivial 0 or 2-curve
and $\c = \b \cup \b_1$ where $\b,\b_1$ are arcs of
$\a,\a_1$ resp.. If $\c$ is a 0-curve then either $D$ or $D^c$
is a trivial disk. If $D$ is trivial then  can isotope to remove
two intersections of $\a$ and $\a_1$ by pulling $\b$ through
$\b_1$. If $D^c$ is trivial then
any of the other 3 disks with boundary consisting of one arc of
$\a$ and $\a_1$ are trivial and hence can reduce the number of
intersections as in first case.

If $\c$ is a 2-curve then each of $\b,\b_1$   intersect
$G$ as if say $\b$ didn't then it would be a compression
for $\a_1$ which contradicts $\a_1$ being non-trivial. Therefore
either $D$ or $D^c$ is trivial i.e. a \nbd of an arc of $G$.
If $D$ is trivial then can isotope by pulling $\b$ through $\b_1$
reducing number of intersections of $\a$ and $\a_1$. If $D^c$ is
trivial then other 3 disks with boundary consisting of one arc
of $\a$ and $\a_1$ cannot have boundaries being n-curves $n \geq 4$
as they can have a maximum of 10 intersections with $G$ between them.
Therefore there is region $\bar{D}$ which either doesn't intersect
$G$ and thus we can isotope as before to reduce the number of 
intersections of $\a$ and $\a_1$ or is \nbd of an arc of $G$ which 
can also be isotoped as before.

So we can assume $\a$ and $\a_1$ intersect twice and divide 
$S^2$ into 4 disks. We label these disks $D_i,i = 1,\ldots,4$ where
$D_1,D_2,D_3$ contain $v_1,v_2,v$ resp. and $\c^i = \partial{D_i}$.
$\c^i$ is an $n_i$-curve where $\sum n_i = 24$ and $n_i \geq 4$.
Also $\c^i = \b^i \cup \b^i_1$ where $\b^i,\b^i_1$ are arcs of
$\a,\a_1$ resp.(figure \ref{setup}). Note that the arcs $\b^i,\b^i_1$
have duplication with each of two arcs that $\a$ or $\a_1$ is divided repeated
twice. This is for ease of labeling and can be thought of as the two sides
of the same arc on $\a$ or $\a_1$.

\begin{figure}
\begin{center}~
\psfig{file=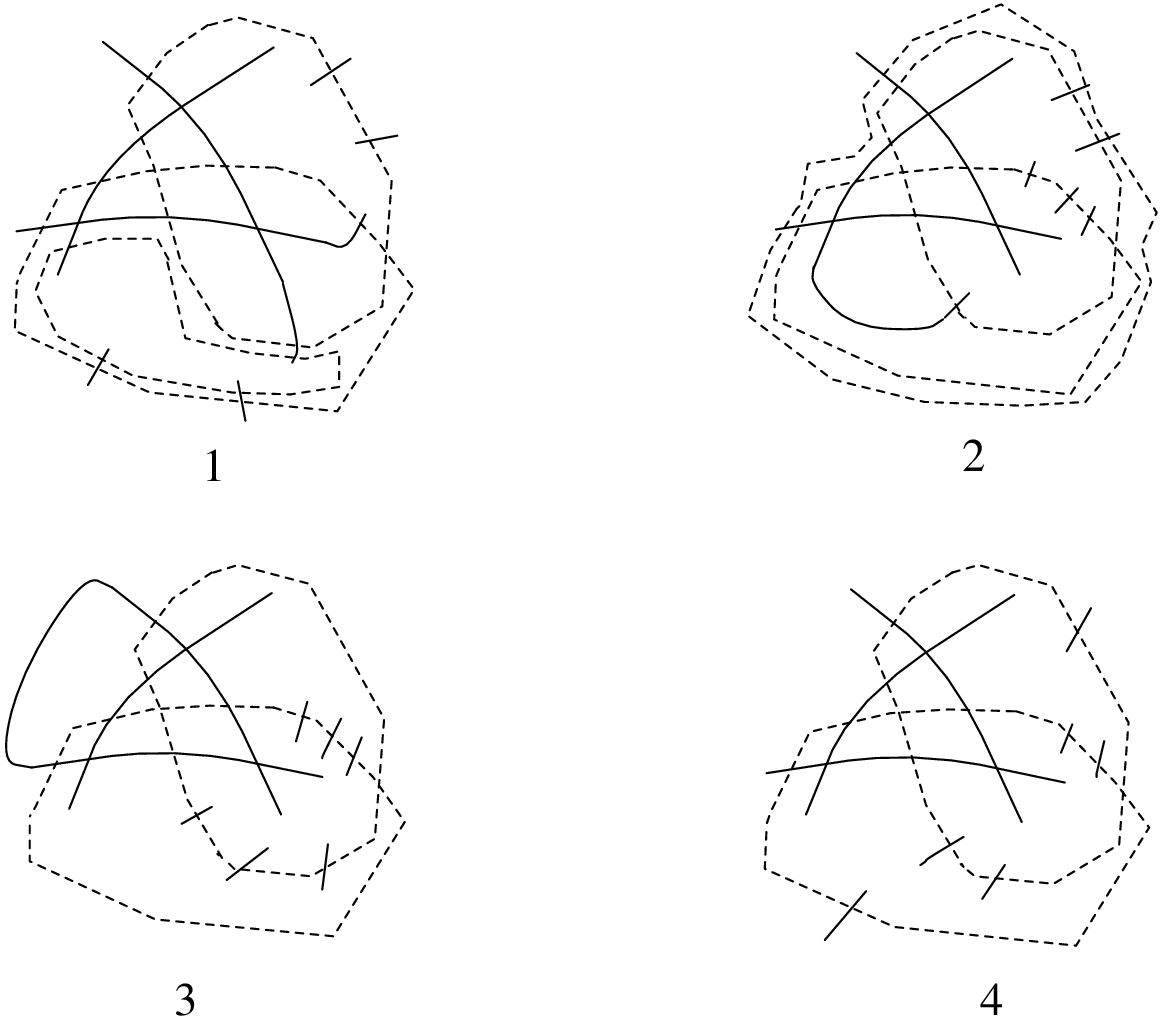,height=2.0in}
\end{center}
\caption{Cases
\label{cases}}
\end{figure}

{\bf Case 1:}If $n_2 = 4$ then  $D_2$ is a \nbd of $v_2$ and both $\b^2$ and
$\b^2_1$  intersect $G$ twice otherwise we get a compression of
$\a$ or $\a_1$. Therefore $\c^1$ is a 6-curve and inside $D_1$ is
a 4-curve $\tilde{\c}^1$(figure \ref{cases}). If it is boundary of
 \nbd of
two parallel arcs of $G$ then this implies either $F_3$ is a bigon
or $\a_1$ is a trivial 6-curve. If it is boundary of
 \nbd of a vertex then
this would imply that $F_3$ was triangular. Therefore $n_2 \neq 4$.

{\bf Case 2:}If $n_1 = 4$ then $\b^1_1$ intersects $G$ only twice as otherwise
$\b^1$ is a compression of $\a_1$. Therefore $\c^4$ is a 6-curve
and $D_4$ contains a 4-curve $\tilde{\c}^4$(figure \ref{cases}).
If it is boundary of \nbd of two parallel arcs of $G$ then this
implies that either $F_1$ is  a bigon or $\a$ is trivial 6-curve.
If it is the boundary of \nbd of vertex then this implies that 
$F_1$ is triangular. Therefore $n_1 \neq 4$ and by symmetry
$n_3 \neq 4$.

{\bf Case 3:}If $n_4 = 4$ then $D_4$ must be \nbd of parallel arcs of $G$
which implies $F_1$ is a bigon(figure  \ref{cases}).
Therefore $n_4 \neq 4$.

{\bf Case 4:}Therefore $n_i = 6$ and each arc $\b^i,\b^i_1$ intersects $G$
exactly 3 times. Therefore $\tilde{\c}^1,\tilde{\c}^4$ are both
4-curves. If $\tilde{\c}^4$ is boundary of \nbd of two parallel 
arcs of $G$ then this implies either $F_1$ is a bigon or both
$\a$ and $\a_1$ have compressions, contradicting them being
non-trivial(figure \ref{cases}). If $\tilde{\c}^4$ is
boundary of a \nbd of a vertex  then $F_1$  would be triangular.
Therefore we have shown that there cannot exist 6-curves 
intersecting as $\a$ and $\a_1$ do  $\Rightarrow$ if $\a$ exists
(i.e. $G'$ isn't atoroidal) $\Rightarrow \a_1$ can't exist 
$\Rightarrow G'_1$ is atoroidal. Similarly $G'_2$ is atoroidal also \eproof

We have shown that the initial objects of $\prec$ are $\{T_n\}_{n \geq 3}$
along with the 4 exceptions. $T_n$ is the projection of the (3,n) torus
link with the link having three components if 3 divides n and having one 
component otherwise. From this we see that $T_n$ has
symmetries taking any directed edge of one of the non-triangular faces to
any other. Therefore any surgery on $T_n$ gives the same graph which we call
$T_n^{+}$.

\begin{figure}
\begin{center}~
\psfig{file=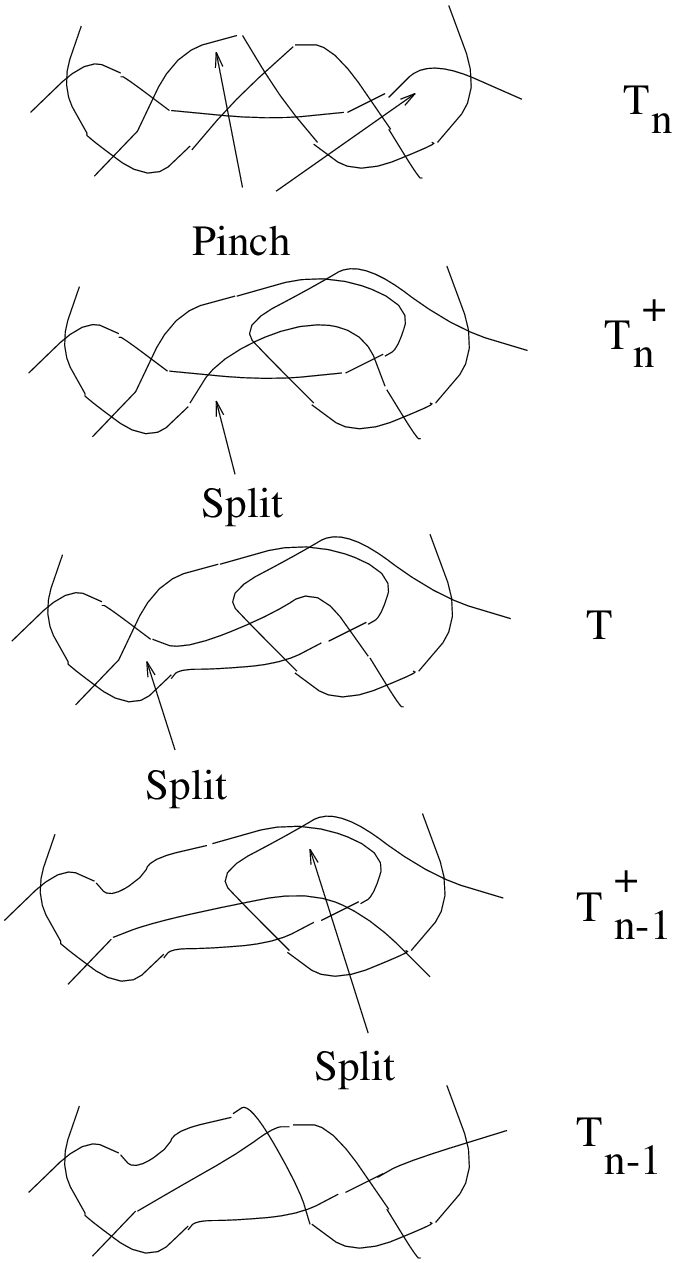,height=2.75in}
\end{center}
\caption{$T^+_{n-1} \prec T^+_{n}$
\label{order}}
\end{figure}

\begin{lemma}
If $T_n \prec G$ and $T_n \neq G$ ($n > 4$)  then  $T_{n-1} \prec G$.
\end{lemma}
\pf If $T_n \prec G$ and $T_n \neq G$ then $T_n^+ \prec G$. 
As in figure \ref{t4}
we can pinch together edges of $T^l_2$ to $T^r_4$ to get $T_n^+$. Only one
vertex of $T_1$ is simple so we cut open at that vertex first. This
reduces the pinched $T_2^l$ to a triangle which has only one simple vertex 
which we now cut open(figure \ref{order}).
This resulting graph is $T_{n-1}^+$ so
we have that $T_{n-1} \prec T_{n-1}^+ \prec T_n^+ \prec G$ \eproof \newline

Therefore if $C_n$ is the set of proper descendants of $T_n$(i.e. 
$T_n \not\in C_n$) then  
$$C_{4} \supseteq C_{5} \supseteq C_{6} \supseteq \cdots \supseteq C_{n}
 \cdots$$ 

\section{Enumeration}
\begin{figure}
\begin{center}~
\psfig{file=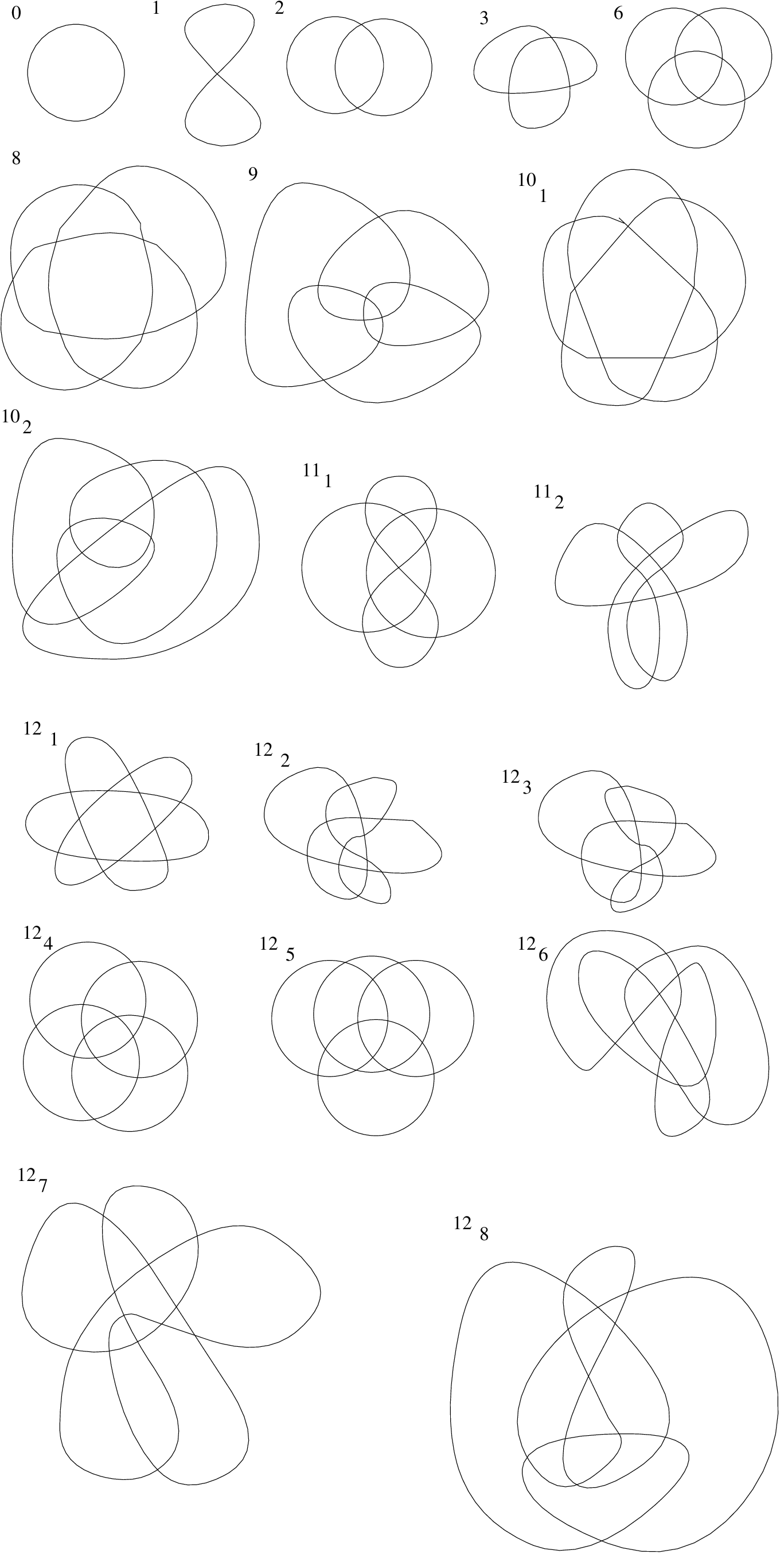,height=6.5in}
\end{center}
\caption{Atoroidal Graphs of 12 crossings or less \label{12crossings}}
\end{figure}

We now have a way to enumerate atoroidal graphs  up to any prescribed crossing
number by performing surgery on the initial objects. Figure \ref{12crossings}
is the enumeration of atoroidal graphs up to 12 crossings. To 
enumerate prime link projections we need only recombine the atoroidal graphs as
follows. We choose two atoroidal graphs $G_1$ and $G_2$ with vertices $v_1$ 
and $v_2$ respectively. Now take the compliment of a neighborhood of 
each vertex and attach their boundaries, making sure to match up the strands 
of the graphs. In recombining we do not 
use the first 3 exceptions as either they have no vertices or the 
compliment of a  neighborhood of a vertex is trivial. 
To enumerate the basic polyhedra of 
Conway the trefoil projection is also not used as the compliment of a vertex is
a bigon.

\end{document}